\documentclass[12pt,leqno]{article}

\usepackage{amsmath,amssymb,amsthm,amsfonts,mathrsfs,amscd,environ}
\usepackage{latexsym,enumerate,color,geometry}
\usepackage{fancyhdr}
\usepackage{verbatim}
\usepackage{hyperref}
\selectfont

 \NewEnviron{ews}{%
\begin{equation}\begin{split}
  \BODY
\end{split}\end{equation}
}

\NewEnviron{ews*}{%
\begin{equation*}\begin{split}
  \BODY
\end{split}\end{equation*}
}

\def\beg{\begin}
\def\be{\begin{equation}}
\def\ee{\end{equation}}
\def\bs{\begin{split}}
\def\es{\end{split}}
\def\bde{\begin{definition}}
\def\ede{\end{definition}}
\def\bews{\begin{ews}}
\def\nws{\end{ews}}

\def\lb{\label}
\def\d{\mathrm{d}}

 \newtheorem{lemma}{Lemma}[section]
 \newtheorem{theorem}[lemma]{Theorem}
 \newtheorem{corollary}[lemma]{Corollary}
 \newtheorem{remark}[lemma]{Remark}
 
 \newtheorem{definition}{Definition}
 \newtheorem{prop}[lemma]{Proposition}
 \numberwithin{equation}{section}

\def\lb{\left(}
\def\rb{\right)}
\def\ra{\rightarrow}
\def\ff{\frac}
\def\<{\langle}
\def\>{\rangle}
\def\mbfR{\mathbf{R}}
\def\vv{\epsilon}
\def\RR{\mathbb{R}}
\def\PP{\mathbb{P}}
\def\E{\mathbb{E}}
\def\Dom{\mathscr{D}}

\title{{\bf Shift Harnack Inequality and Integration by Part Formula for Semilinear SPDE}\footnote{Supported by NSFC(11131003), SRFDP, 985-Project.}}

\author{
{\bf Shao-Qin Zhang }\\
\footnotesize{School of Math. Sci. and Lab. Math. Com. Sys., Beijing Normal University, Beijing 100875, China}\\
\footnotesize{Email: zhangsq@mail.bnu.edu.cn}
}

\begin{document}
\maketitle

\begin{abstract}
Shift Harnack and integration by part formula are establish for semilinear spde with delay and a class of stochastic semilinear evolution equation which cover the hyperdissipative Naiver-Stokes/Burges equation. For the case of stochastic equation with delay, an extension to path space is given.
\end{abstract}\noindent

AMS Subject Classification:\ 60H15
\noindent

Keywords: Shift Harnack inequality, integration by part formula, stochastic functional equations, hyperdissipative Navier-Stokes/Burges equation, path spaces, log Sobolev inequality.

\vskip 2cm

\section{Introduction}
Recently, using a new coupling argument, \cite{Wang2012} provides a new type Harnack inequality, called shift Harnack inequality, and derive Driver's integration by part formula, see \cite{Dri1997}. The main idea is that construct two processes which start from the same point and at the expected time T they separate at a fixed vector almost surely. In \cite{Wang2012} there, for the case of semilinear stochastic partial differential equations(SPDE), two problems remains, the first one is that how to establish shift Harnack inequality and integration by part formula for semilinear SPDE with delay, the second is that whether the two processes can separate at arbitrarily vector. In this paper, we try to find the answer to the two problems. We construct coupling in the spirit of \cite{Wang2012}, for the case of stochastic functional equation it even dates back to \cite{ERS2009}. With a little knowledge of control theory and regularity theory of semigroups, explicit coupling is constructed, then we derive the shift Harnack and integration by part formula.

In the second part of the paper, we deal with semilinear SPDE with delay and generalized to non-Lipschitz case. In the third part, we extend the integration by part formula to the path space of the solution of the stochastic functional equation, some application are given. The last part, we establish the shift Harnack inequality and integration by part formula for a class of stochastic evolution equation, which covers the hyperdissipative Navier-Stokes/ Burgers equation.

Before our main results, we need some preparation. For any $T>0$, assume that $U,H$ are Hilbert spaces, $\ B\in \mathscr{L}(U,H)$, and $-A$ generates an analytic semigroup, define two operators as follow, which are well know in control theory,
\begin{eqnarray}
&&L^B_T:L^{2}([0,T],U)\ra H,\ L^B_Tf=\int_0^T{e^{-(T-t)A}Bf(t)\d t},\\
&&\mathbf{R}^B_T:H\ra H,\ \mathbf{R}^B_Th=\int_0^T{e^{-tA}BB^*e^{-tA^*}h\d t},
\end{eqnarray}
and let
\be
D_A(1/2,2)=\{x\in H\ |\ ||x||^2_\frac{1}{2}:=\int_0^\infty||Ae^{-tA}x||^2\d t<\infty\}.
\ee

The following proposition are well know in semigroup theory and control theory, see Theorem 3.1 in page 143 in \cite{BDPDM2007} and Appendix B in \cite{DPZ1992} for details,
\begin{prop}\label{prop1}
(1) Assume that $-A$ generate an analytic semigroup, then for each $T>0$, the map
\be
u\ra (u'+Au,u(0)): W^{1,2}([0,T],H)\bigcap L^2([0,T],\Dom(A))\ra L^{2}([0,T],H)\times D_A(1/2,2 ),
\ee
is an isomorphism, and
\be
W^{1,2}([0,T],H)\bigcap L^2([0,T],\Dom(A))\subset C([0,T],D_A( 1/2,2 )).
\ee
(2) For the two operator $L^B_T$ and $\mbfR^B_T$, if $B^{-1}\in \mathscr{L}(H,U)$, then
\begin{eqnarray}
\textbf{Im}(L^B_T)&=&\textbf{Im}((\mbfR^B_T)^{\frac{1}{2}})=D_A(1/2,2),\\
||(\mbfR^B_T)^{-1/2}x||^2&=& \min{\{\int_0^T||f(s)||_U^2\d s\ |\ L_Tf=x\}},\  x\in D_A(1/2,2).
\end{eqnarray}
Here  $(\mbfR^B_T)^{-1/2}$ means the pseudo-inverse of $(\mbfR^B_T)^{1/2}$.
\end{prop}
\begin{remark}\label{lmm1}
By Proposition \ref{prop1}, for any $T>0$, $x\in D_A(1/2,2)$, there exists $f\in L^{2}([0,T],U)$ such that
\be
L^B_Tf=x,\ ||f||^2_{L^2([0,T],U)}=\min{\{\int_0^T||f(s)||_U^2\d s\ |\ L^B_Tf=x\}}.
\ee
\end{remark}
The following lemma will give us the time behavior
\begin{lemma}\label{lmm3}
If $B^{-1}\in \mathscr{L}(H,U)$ and $Ae^{-tA}(D_A(1/2,2))\subset U$ for all $t>0$, then exists $C>0$ which independent of T, such that
\be
||(\mbfR^B_T)^{-1/2}x||^2\leq 2||B^{-1}||^2\Big(\ff {||x||^2} {T} + 2||x||^2_{\ff 1 2}\Big)
\ee
\end{lemma}
\beg{proof} For all $x\in D_A(1/2,2)$,
let $f(t)=B^{-1}\Big(\ff {e^{-(T-t)A}x} {T} + \ff {2t} {T}Ae^{-(T-t)A}x\Big)$. Then
$L_Tf=x$ and
\beg{ews}
&||(\mbfR^B_T)^{-1/2}x||^2\leq ||f||^2_{L^2([0,T],U)}
\leq ||B^{-1}||^2\Big(\ff {2||x||^2} {T} + 4 \int_0^T||Ae^{-(T-t)A}x||^2\d t\Big)\\
&\leq 2||B^{-1}||^2\Big(\ff {||x||^2} {T} + 2||x||^2_{\ff 1 2}\Big).
\end{ews}
\end{proof}

\section{Semilinear SPDE with Delay}
$H$ is Hilbert space with norm $||\cdot||$, $\mathscr{C}=C([-\tau,0],H)$, consider the following equation 
\be\label{equ0}
\d x(t)=-Ax(t)\d t +F(x_t)\d t +\sigma(t)\d W(t),
\ee
satisfies the following conditions
\begin{enumerate}[({H}1)]
\item $-A$ generates analytic semigroup, there exists $a\in \RR$ such that $-A-a$ is dissipative, and there exists $\alpha\in(0,\frac{1}{2})$ so that $\int_0^Tt^{-2\alpha}||e^{-tA}||^2_{HS}\d t<\infty$,\label{H1}
\item $F:\mathscr{C}\ra H$ is Lipschitz with Lipschitz constant $L$,\label{H2}
\item $\sigma:[0,T]\ra \mathscr{L}(H)$ measurable and bounded, and there is $M>0$, such that $||\sigma(\cdot)^{-1}||_\infty\leq M$.\label{H3}
\end{enumerate}
Denote the solution of the equation with initial value $\xi$ by $x(t,\xi)$, related segment process $x_t(\xi)$, and $P_T f(\xi)=\E f(x_T(\xi))$. In this section, we choose that $B=I$ in Proposition \ref{prop1}. Now, we state our result
\begin{theorem}\label{theorem1}
Assume that (H\ref{H1}) to (H\ref{H3}) hold and $T>\tau$. Let $$\eta\in W^{1,2}([-\tau,0],H)\bigcap L^2([-\tau,0],\Dom(A)),\ \psi(t)=\eta'(t-T)+A\eta(t-T),\ t\in[T-\tau,T].$$ Then
for any $\xi\in \mathscr{C}$, $f\in\mathscr{B}_b(\mathscr{C})$,
\begin{ews}\label{shift_H_1}
\lb P_T f(\xi)\rb^p\leq & P_T f^p(\cdot+\eta)(\xi)
\exp\Big\{\frac{pM^2}{p-1}\Big[\frac{2+L^2(T-\tau)^2}{2}||(\mbfR^I_{T-\tau})^{-1/2}\eta(-\tau)||^2\\
&+\int_{T-\tau}^T{||\psi(t)||^2}\d t  +\tau\Big\{\Big[(T-\tau)||(\mbfR^I_{T-\tau})^{-1/2}\eta(-\tau)||^2\Big]\vee||\eta||_\infty^2\Big\}
\Big]\Big\}
\end{ews}
Moreover, if we assume that $F:\mathscr{C}\ra H$ is G\^{a}teaux derivable with $||\nabla F(\cdot)||_\infty\leq L<\infty$ in addition, then for any $\phi\in L^{2}([0,T-\tau],H)$ such that $L^I_{T-\tau}=\eta(-\tau)$, and $\psi$ defined as above, we have
\begin{ews}
&(P_T\nabla_\eta f)(\xi)\\
=&\E \left\{f(x_T(\xi))\int_0^T\<\sigma(t)^{-1}\left[\phi(t)1_{[0,T-\tau)}(t)
+\psi(t)1_{[T-\tau,T]}(t)-\nabla_{\Gamma_t}F(x_t(\xi))\right],\d W(t)\>\right\},\ f\in C_b^1(\mathscr{C}),
\end{ews}
where
\be
\Gamma(t)=\left\{
\begin{array}{lc}
\eta(t-T), & t\geq T-\tau,\\
\int_0^te^{-(t-s)A}\phi(s)\d s, & t<T-\tau.\\
\end{array}\right.
\ee
\end{theorem}

\noindent\emph{Proof of Theorem \ref{theorem1}.} By Proposition \ref{prop1} for the case that $B=I$ and $U=H$, $\psi$ are well defined. By Remark \ref{lmm1}, firstly we choose $\phi\in L^{2}([0,T-\tau],H) $ such that
\begin{ews}\label{equ6}
L^I_{T-\tau}\phi=\eta(-\tau),\ \  ||\phi||_{L^{2}([0,T-\tau],H)}=||(\mbfR^I_{T-\tau})^{-1/2}\eta(-\tau)||.
\end{ews}
We construct another process as follow
\begin{eqnarray*}
\left\{
\begin{array}{l}
\d y(t)=-Ay(t)\d t+F(x_t(\xi))\d t +\sigma(t)\d W(t) + \epsilon (\phi(t)1_{[0,T-\tau)}(t)\d t+ \psi(t)1_{[T-\tau,T]}(t))\d t,\\
y_0=\xi,\\
\end{array}\right.
\end{eqnarray*}
then
\be
y(t)=x(t)+\epsilon\int_0^t{e^{-(t-s)A}\phi(s)1_{[0,T-\tau)}(s)}\d s+\epsilon\int_0^t{e^{-(t-s)A}\psi(s)1_{[T-\tau,T]}(s)}\d s.
\ee
For $t\geq T-\tau$,
\begin{ews}
\int_0^t{e^{-(t-s)A}\phi(s)1_{[0,T-\tau)}(s)}\d s =&e^{-(t-T+\tau)A}\int_0^{T-\tau}{e^{-(T-\tau-s)A}\phi(s)}\d s\\
=&e^{-(t-T+\tau)A}\eta(-\tau).
\end{ews}
Since $\eta\in W^{1,2}([-\tau,0],H)\bigcap L^2([-\tau,0],\Dom(A))$, and
\be
\psi(t)=\eta'(t-T)+A\eta(t-T),\ t\geq T-\tau,
\ee
that means $\eta(\cdot-T)$ is the solution of the following equation
\be
\frac{\d \eta(t-T)}{\d t}=-A\eta(t-T)+\psi(t), t\geq T-\tau,
\ee
with initial value $\eta(-\tau)$ at $T-\tau$, or in the integration form
\be
\eta(t-T)=e^{-(t-T+\tau)}\eta(-\tau)+\int_{T-\tau}^{t}e^{-(t-s)A}\psi(s)\d s,\ t\geq T-\tau,
\ee
thus, for $t\geq T-\tau$,
\begin{ews}\label{equ5}
y(t)=&x(t)+\epsilon e^{-(t-T+\tau)A}\eta(-\tau)+ \epsilon(\eta(t-T)-e^{-(t-T+\tau)}\eta(-\tau))\\
=&x(t)+\epsilon\eta(t-T),
\end{ews}
that means
\be
y_T=x_T+\epsilon\eta,
\ee
therefore, for all $t\in [0,T]$,
\begin{eqnarray}
y(t)-x(t)=\epsilon\Gamma(t),\ y_t-x_t=\epsilon\Gamma_t.
\end{eqnarray}
Let
$$h^\vv(t)=\epsilon\sigma(t)^{-1}\lb\phi(t)1_{[0,T-\tau)}(t)+\psi(t)1_{[T-\tau,T]}(t)\rb +\sigma(t)^{-1}(F(x_t)-F(y_t)),$$
$$\d \tilde{W}(t)=\d W(t)+h^\epsilon(t)\d t,\ \ R^\epsilon_T=\exp\left\{-\int_0^T\<h^\vv(t),\d W(t)\>-\frac{1}{2}\int_0^T||h^\epsilon(t)||^2\d t\right\}.$$
Then we can rewrite the equation of $y$ as
\be
\d y(t)=-Ay(t)\d t + F(y_t)\d t+ \sigma(t)\d \tilde{W}(t),\ y_0=\xi.
\ee
By (H\ref{H1}) to (H\ref{H3}), as in \cite{Wang2012}, and noting that
\begin{ews}
\int_0^T||\Gamma_t||_\infty^2\d t &\leq \int_0^{T-{\tau}}\lb\int_0^t||\phi(s)||\d s\rb^2\d t+\int_{T-\tau}^{T}||\Gamma_t||^2\d t\\
&\leq  \int_0^{T-\tau}t\int_0^{T-\tau}||\phi(s)||^2\d s\d t+\tau\left\{\left[(T-\tau)\int_0^{T-\tau}||\phi(s)||^2\d s\right]\vee||\eta||_\infty^2\right\}\\
&\leq \frac{(T-\tau)^2}{2}||(\mbfR^I_{T-\tau})^{-1/2}\eta(-\tau)||^2
+\tau\left\{\left[(T-\tau)||(\mbfR^I_{T-\tau})^{-1/2}\eta(-\tau)||^2\right]\vee||\eta||_\infty^2\right\},
\end{ews}
we can prove that $\{\tilde{W}(t)\}_{t\in[0,T]}$ is Brownian motion by Girsanov theorem and get the shift Harnack inequality
\begin{ews}
\lb P_T f(\xi)\rb^p\leq & P_T f^p(\cdot+\eta)(\xi)(\E(R_T^1)^{\frac{p}{p-1}})^{p-1}\\
\leq & P_T f^p(\cdot+\eta)(\xi)\exp\Big[\frac{p}{2(p-1)}\int_0^T||h(t)||^2\d t\Big]\\
\leq & P_T f^p(\cdot+\eta)(\xi)
\exp\Big\{\frac{pM^2}{p-1}\Big[\frac{2+L^2(T-\tau)^2}{2}||(\mbfR^I_{T-\tau})^{-1/2}\eta(-\tau)||^2\\
&+\int_{T-\tau}^T{||\psi(t)||^2}\d t  +\tau\Big\{\Big[(T-\tau)||(\mbfR^I_{T-\tau})^{-1/2}\eta(-\tau)||^2\Big]\vee||\eta||_\infty^2\Big\}
\Big]\Big\}
\end{ews}
in the last inequality, we have used that (\ref{equ6}).
Further more since $F$ has bounded G\^{a}teaux derivative, choosing $\phi\in L^{2}([0,T-\tau],H) $ such that
$L^I_{T-\tau}\phi=\eta(-\tau)$, then
\begin{ews}
\frac{\d}{\d \epsilon}|_{\epsilon=0}R_T^\epsilon =&-\int_0^T\<\sigma(t)^{-1}\left[\phi(t)1_{[0,T-\tau)}(t)
+\psi(t)1_{[T-\tau,T]}(t)\right],\d W(t)\>\\
&+\int_0^T\<\sigma(t)^{-1}\nabla_{\Gamma_t}F(x_t),\d W(t)\>,
\end{ews}
hold in $L^1(\PP)$. Therefore
\begin{ews}
&(P_T\nabla_\eta f)(\xi)=\lim_{\vv\ra 0^+}\Big[\frac{P_Tf(\cdot-\vv\eta)(\xi)-\E f(x_T(\xi))}{-\vv}\Big]\\
&= \lim_{\vv\ra 0^+}\Big[\frac{\E R^\vv f(y_T(\xi)-\vv\eta)-\E f(x_T(\xi))}{-\vv}\Big]= \lim_{\vv\ra 0^+}\Big[\frac{\E R^\vv f(x_T(\xi))-\E f(x_T(\xi))}{-\vv}\Big]\\
&=\E \left\{f(x_T(\xi))\int_0^T\<\sigma(t)^{-1}\left[\phi(t)1_{[0,T-\tau)}(t)
+\psi(t)1_{[T-\tau,T]}(t)-\nabla_{\Gamma_t}F(x_t(\xi))\right],\d W(t)\>\right\}.
\end{ews}
\qed

\beg{remark}
The second condition in (\ref{equ6}) is only used to shift Harnack inequality, to get explicit integration by part formula, one can choose $``\phi "$, here we provide a procedure to get it and give an example. Fix any $T>0$. For any $x\in D_A(1/2)$, $h\in L^2([0,T],H)$, let
$$\phi_1(t)=e^{-tA}x+\int_0^t e^{-(t-s)A}h(s)\d s,$$ Then  $\phi_1\in W^{1,2}([0,T],H)\bigcap L^2([0,T],\Dom(A)$ by Proposition \ref{prop1}. Let $u\in C^1([0,T],\mathbb{R})$, $u(0)=0,\ u(T)=1$. Then $$\phi(t)=\ff {\d} {\d t}\Big(u(t)\phi_1(T-t)\Big)+u(t)A\phi_1(T-t).$$
It's clear that $\phi\in L^2([0,T],H)$ and
$$\int_0^T e^{-(T-t)A}\phi(t)\d t=u(T)\phi_1(0)-u(0)\phi_1(T)=x.$$
For example, one can choose $h=0$, $u(t)=\ff t T$, then $\phi(t)=\ff {e^{-(T-t)A}x} {T} + \ff {2t} {T}Ae^{-(T-t)A}x$
\end{remark}

For general case we can use Lemma \ref{lmm3} to get the following inequality, for more sharp estimate we expect more better inequality.
\beg{corollary}
Assume that (H\ref{H1}) to (H\ref{H3}) hold, $T>\tau$ and $\psi$ as in Theorem \ref{theorem1}, then
\begin{ews}\label{shift_H_1}
\lb P_T f(\xi)\rb^p\leq & P_T f^p(\cdot+\eta)(\xi)
\exp\Big\{\frac{pM^2}{p-1}\Big[\Big(2+L^2(T-\tau)^2\Big)\Big(\ff {||x||^2} {T-\tau} + 2||x||^2_{\ff 1 2}\Big)\\
&+\int_{T-\tau}^T{||\psi(t)||^2}\d t  +\tau\Big\{\Big[2(T-\tau)\Big(\ff {||x||^2} {T-\tau} + 2||x||^2_{\ff 1 2}\Big)\Big]\vee||\eta||_\infty^2\Big\}
\Big]\Big\}
\end{ews}
\end{corollary}

When $A$ is self adjoint, we have
\begin{corollary}
When $A$ is self adjoint operator and $A\geq \lambda_0>0$, under the assumption in Theorem \ref{theorem1}, we have the shift Harnack inequality holds with $||(\mbfR^I_{T-\tau})^{-1/2}||^2$ replaced by $\ff {2||A^{1/2}\eta(-\tau)||^2} {1-e^{-2(T-\tau)\lambda_0}}$.
\end{corollary}
\noindent\emph{Proof.}
In this situation,
\be
\mbfR^I_{T-\tau}=\int_0^{T-\tau} e^{-2tA}\d t= \frac{A^{-1}}{2}(I-e^{-2(T-\tau)A}),
\ee
and $D_A(1/2,2)=\Dom(A^{1/2})$, then
\be
||(\mbfR^I_{T-\tau})^{-1/2}\eta(-\tau)||^2=2||(I-e^{-2(T-\tau)A})^{-\frac{1}{2}}A^{\frac{1}{2}}\eta(-\tau)||^2\leq \frac{2\left|\left|A^{1/2}\eta(-\tau)\right|\right|^2}{1-e^{-2(T-\tau)\lambda_0}}.
\ee
\qed
\begin{corollary}\label{coro3}
Assume that (H\ref{H1}), (H\ref{H3}) hold and there is an increasing continuous function  $\gamma:[0,\infty)\ra [0,\infty)$ such that                
\be
||F(x)-F(y)||\leq \gamma(||x-y||_\infty),\ \forall x,y\in\mathscr{C},
\ee
and equation (\ref{equ0}) has pathwise unique mild solution, then for any $p>1$, $f\in \mathscr{B}_b^+(H)$,
\begin{ews}
(P_Tf)^p &\leq P_Tf^p(e+\cdot) \exp\Big[\frac{M^2p}{p-1}\Big(||(\mbfR^I_{T-\tau})^{-1/2}\eta(-\tau)||^2\\
&\qquad+T\gamma^2\lb||\eta||_\infty\vee\sqrt{(T-\tau)||(\mbfR^I_{T-\tau})^{-1/2}\eta(-\tau)||^2}\rb+\int_{T-\tau}^T||\psi(t)||^2\d t\Big)\Big].
\end{ews}
\end{corollary}
\noindent\emph{Proof.}
We use the notation in Theorem \ref{theorem1}. In this case, the shift Harnack inequality follows from the following estimate
\begin{ews}
\frac{1}{2M^2}\int_0^T||h(t)||^2\d t
\leq & \lb\int_0^{T-\tau}||\phi(t)||^2\d t+\int_{T-\tau}^T||\psi(t)||^2\d t+\int_0^T\gamma^2(||\Gamma_t||)\d t\rb \\
\leq & ||(\mbfR^I_{T-\tau})^{-1/2}\eta(-\tau)||^2+\int_{T-\tau}^T||\psi(t)||^2\d t\\
& +\int_0^T\gamma^2\lb\int_0^{t\wedge(T-\tau)}||\phi(s)||\d s\vee||\eta||_\infty\rb\d t\\
\leq & ||(\mbfR^I_{T-\tau})^{-1/2}\eta(-\tau)||^2
+\int_{T-\tau}^T||\psi(t)||^2\d t\\
&+T\gamma^2\lb||\eta||_\infty\vee\sqrt{(T-\tau)||(\mbfR^I_{T-\tau})^{-1/2}\eta(-\tau)||^2}\rb.
\end{ews}
\qed

\begin{remark}
For the existence and uniqueness of stochastic functional differential equations with non-Lipschitz coefficients and nontrivial examples, one can see \cite{BH2010,ShaoWY2012,WangZTS2010} and references there in.
\end{remark}

\section{Extend to Path Space}

Firstly, we shall give an integration by part formula on path space of solution of stochastic functional differential equation, it follows from \cite{Hsu2002, Wang11}. For simplicity, we assume that $\nabla F$ is bounded, $\sigma\equiv I$. Denote $\nabla^j$ the partial derivative of the j-th component, $W_T^\xi$  the path space of segment process $x_\cdot^\xi$ on $[0,T]$. Denote all the bounded smooth cylindrical function on $\mathscr{C}$ by $\mathscr{F}C_b^\infty(W^\xi(T))$, i.e.
$$\mathscr{F}C_b^\infty(\mathscr{C})=\big\{G(\gamma)=g(\gamma_{t_1},\cdots,\gamma_{t_n}),\ g\ \mbox{is smooth on}\ \mathscr{C}^n\ \mbox{with bounded all derivative },\ \gamma\in W^\xi(T)  \big\}$$

\beg{definition}
Let G be a function on $W^\xi(T)$, it's directive along the direction $\eta$ at $\gamma\in W^\xi(T)$, if the following limit exists $$\nabla_\eta G(\gamma):= \lim_{\vv\ra 0^+}\ff {G(\gamma+\vv\eta)-G(\gamma)} {\vv}.$$ If $\nabla_\cdot G(\gamma)$ provides a continuous linear functional on $W^{1,2}_0([0,T],H):=\Big\{f\in W^([0,T],H)| f(0)=0\Big\}$, the gradient $\nabla G (\gamma)$ is defined as its Riesz representation, i.e.  $\nabla G(\gamma)$ is an element in $W^{1,2}_0([0,T],H)$ such that $\<\nabla G(\gamma), \eta\>_{W^{1,2}_0}=\nabla_\eta G(\gamma)$.
\end{definition}

By definition, for $G(\gamma)=g(\gamma(t_1), \cdots,\gamma(t_n))$, $g\in C_b^1(H^n)$, it's easy to see that $\ff {\d} {\d t}(\nabla G(\gamma))(t)=\sum_{i=1}^{n}1_{[t<t_i]}\nabla^i g(\gamma)$.

\beg{prop}
If $F$ is Fr\'ech\`et differentiable on $\mathscr{C}$ with $||\nabla F||\leq L$ and $\nabla F$ is uniformly continuous on bounded set of $\mathscr{C}\times\mathscr{C}$. Then for all $\eta\in W^{1,2}_0([0,T],H)\bigcap L^2([0,T],\Dom(A)) $, extending it to $[-\tau, 0]$ by zero, we have
$$\E \nabla_\eta G(x^\xi([0,T]))=\E G(x^\xi([0,T]))\int_0^T\<\ff {\d \eta(t)}{\d t}+A\eta(t)-\nabla_{\eta_t}F(x_t),\d W(t)\>,\ G\in C_b^1(W^\xi(T)).$$
\end{prop}
\beg{proof}
By the definition of $\eta$, we can definite
$$h(t)=\eta(t)+\int_0^t A\eta(s)\d s-\int_0^t\nabla_{\eta_s}F(x_s)\d s.$$ Then $h\in W^{1,2}_0([0,T],H)$ $\PP$-a.s.  and it's adapted. Since $\nabla F$ is Fr\'ech\`et differentiable on $\mathscr{C}$ with $||\nabla F||\leq L$ and $\nabla F$ is uniformly continuous on bounded set of $\mathscr{C}\times\mathscr{C}$,  Let $D_h x$ be the Malliavian derivative of $x$ along $h$, then is the mild solution of the following equation
$$\d D_h x(t)=-AD_h x(t)\d t+ \nabla_{D_h x_t}F(x_t)\d t+h'(t)\d t,\ D_h x_0=0,$$
by \cite[Theorem A.2]{BYWang}. Noting that $\eta$ and $D_h x$ satisfy the same differential equation with the same initial value, then $D_h x=\eta$ $\PP$-a.s. Therefore
\beg{ews}
\E\nabla_\eta G(x^\xi([0,T]))&=\E\nabla_{D_h x}G(x^\xi([0,T]))=\E D_h G(x^\xi([0,T]))\\
&=\E G(x^\xi([0,T]))\int_0^T\<h(t),\d W(t)\>\\
&=\E G(x^\xi([0,T]))\int_0^T\<\ff {\d \eta} {\d t}+A\eta(t)-\nabla_{\eta_t}F(x_t),\d W(t)\>.
\end{ews}
\end{proof}

Let $\mathscr{E}^\xi(\Phi,\Psi)=\E \<\Phi,\Psi\>_{W^{1,2}_0}(x^\xi([0,T]))$, $\Phi,\Psi\in \mathscr{F}C_b^\infty(W^\xi(T))$. Then we have
\beg{corollary}
$\Big(\mathscr{E}^\xi,\mathscr{F}C_b^\infty(W^\xi(T))\Big)$ is closable in $L^2(W^\xi(T),\Pi^\xi(T))$.
\end{corollary}
\beg{proof}
Let $\phi \in C_b^2(\mathbb{R}^{n\times m})$ and $$\Phi(\gamma)=\phi(\<\gamma(t_1),e_1\>,\cdots,\<\gamma(t_1),e_n\>,
\cdots,\<\gamma(t_n),e_1\>,\cdots,\<\gamma(t_m),e_n\>),$$
for any $\Psi\in \mathscr{F}C_b^\infty(W^\xi(T))$, by integration by part formula,
\beg{ews}
\E\<\nabla\Phi,\nabla \Psi\>_{W^{1,2}_0}&=\sum_{i,j=1}^{m,n}\E(\partial \phi)_{ij}\<\Psi,\cdot\wedge t_i e_j\>_{W^{1,2}_0}\\&=\sum_{i,j=1}^{n,m}\E\Big(\nabla_{\cdot\wedge t_i e_j}[\Psi(\partial \phi)_{ij}]-\Psi\nabla_{\cdot\wedge t_i e_j}(\partial \phi)_{ij}\Big)\\
&=\E \Psi\sum_{i,j=1}^{m,n}\Big[(\partial \phi)_{ij}\delta(\cdot\wedge t_i e_j)-\nabla_{\cdot\wedge t_i e_j}(\partial \phi)_{ij}\Big].
\end{ews}
By this formula and that this type of $\nabla\Phi$ is dense in $L^2(W^\xi(T)\ra W^{1,2}_0([0,T],H),\Pi^\xi(T))$, we can prove the corollary.
\end{proof}

Next we shall prove the log Sobolev inequality for this Dirichlet form. Let $y(\cdot, h)$ be the the mild solution of the following equation in pathwise sense
\be\label{equ12}\d y(t)=-Ay(t)\d t+\nabla_{y_t}F(x_t)\d t + h'(t)\d t,\ y_0=0,\ee
where $h\in L^2(\Omega\ra W^{1,2}_0([0,T],H))$, then $y \in C([0,T],H)$ $\PP$-a.s., if we assume that $\nabla F$ is bounded, then it's clear that $\nabla_{y_\cdot}F(x_\cdot)+h'(\cdot)\in L^2([0,T],H),$ $\PP$-a.s., by Proposition \ref{prop1}, we find that $y(\cdot, h)\in W^{1,2}_0([0,T],H)\bigcap L^2([0,T],\Dom(A))$, that means it's a continuous operator on $W^{1,2}_0([0,T],H)$. Before the estimate of the operator norm, we recall a priori estimate in \cite[Lemma 3.3 in page 141]{BDPDM2007}.

\beg{lemma}\label{lmm_priori}
Assume that $J$ generates an analytic semigroup with negative type and the resolvent satisfies $||(\lambda-J)^{-1}||\leq M/|\lambda|,\ \textrm{Re}\lambda>0$. For all $f\in W^{1,2}([0,\infty), H)$ and $T>0$. Let $v$ be the solution of the following equation
$$\d v(t)=Jv(t)\d t+f(t)\d t,\ v(0)=0.$$ Then
$$
\int_0^T||Jv(t)||^2\d t \leq (M+1)^2\int_0^T ||f(t)||^2\d t,\ \ \  \int_0^T||v'(t)||^2\d t \leq M^2\int_0^T ||f(t)||^2\d t.
$$
\end{lemma}
\beg{proof}
The first inequality is the result of \cite[Lemma 3.3 in page 141]{BDPDM2007}. The second one was missing there, but it can be proved follow \cite[Lemma 3.3 in page 141]{BDPDM2007} completely, a proof is given here just for convenient.  Let $\bar{v}(t)$ be the solution of
\be\label{equ11} \d \bar{v}(t)=J\bar{v}(t)\d t + \bar{f}(t)\d t,\ \bar{v}(0)=0,\ t\in \RR,\ee
where $\bar{f}(t)=f(t)1_{(0,T)}(t)$. Then
\be
\bar{v}(t)=\left\{\begin{array}{cc}
  v(t) & t\in[0,T] \\
  0 & t\leq 0 \\
  e^{(t-T)J}v(T) & t\geq T.
\end{array}\right.
\ee
Apply Fourier transform to equation (\ref{equ11}), letting
$$\hat{v}(k)=\int_0^\infty e^{-ikt}\bar{v}'(t)\d t,\ \ \hat{f}(k)=\int_0^\infty e^{-ikt}\bar{f}(t)\d t,$$
we arrive at
$$ik\hat{v}(k)=A\hat{v}(k)+\hat{f}(k),$$
then
$$\hat{v}(k)=(ik-J)^{-1}\hat{f}(k),$$
this implies that $$||ik\hat{v}(k)||=||(ik-J)^{-1}\hat{f}(k)+\hat{f}||=||ik(ik-J)^{-1}\hat{f}(k)||\leq M||\hat{f}(k)||,$$
by Parseval's inequality, we get the second inequality.
\end{proof}

\beg{lemma}
Assume that $\sigma\equiv I$, $||\nabla F||\leq L$. Then $y$ is a continuous operator on $W^{1,2}([0,T],H)$ and $$||y(\cdot, h)||^2_{W^{1,2}}\leq e^{2Ta^+ + T^2(a^+)^2e^{2Ta^+}}\Big(1+LTe^{T(L+a^+)}\Big)||h||_{W^{1,2}}^2,\ \ \forall h\in W^{1,2}([0,T],H).$$
\end{lemma}
\beg{proof} It's no harm to assume that $a\geq 0$. Since $y(\cdot,h)\in W^{1,2}_0([0,T],H)\bigcap L^2([0,T],\Dom(A))$, we view $\nabla_{y_\cdot}F(x_\cdot)+h'(\cdot)$ as a inhomogeneous term in equation (\ref{equ12}). Replacing $y(t)$ by $e^{at}(t)$ and $-A$ by $-A-a$, we get
$$\int_0^t||(e^{ra}y(r))'||^2\d r\leq \int_0^te^{2ra}||\nabla_{y_r}F(x_r)+h'(r)||^2\d r,\ \ t\leq T,$$ by Lemma \ref{lmm_priori}. Then
\beg{ews}
&\int_0^t e^{2ra}||y'(r)||^2\d r-\int_0^te^{2ra}||\nabla_{y_r}F(x_r)+h'(r)||^2\d r\\&\leq -\int_0^t||ae^{ra}y(r)||^2\d r-2\int_0^t\<ae^{ra}y(r),e^{ra}y'(r)\>\d r\\
&\leq \int_0^t a^2e^{2ra}||y(r)||^2\d r\leq Ta^2e^{2Ta}\int_0^t\int_0^r||y'(s)||^2\d s\d r.
\end{ews}
By Gronwall's inequality
$$\int_0^T||y'(t)||^2\d t\leq e^{T^2a^2e^{2Ta}+2Ta}\int_0^T(L||y_t||_\infty+||h'||(t))^2\d t.$$
On the other hand
$$\d ||y(t)||^2\leq 2a||y(t)||^2\d t + 2L ||y_t||_\infty||y(t)||\d t+2||h'(t)||\cdot||y(t)||\d t,$$ then
$$\d ||y(t)||\leq a||y(t)||^2\d t + L ||y_t||_\infty\d t+||h'(t)||\d t.$$ This implies that
$$||y_t||_\infty\leq (L+a)\int_0^t||y_r||_\infty\d r + \int_0^t||h'(r)||\d r,\ t\geq 0.$$ By Gronwall's lemma,
$$||y_t||_\infty\leq e^{t(L+a)}\int_0^t||h'(r)||\d r.$$ By this estimate, one can find that
$$L^2 \int_0^T||y_t||_\infty^2 \leq L^2\int_0^Te^{t(L+a)}(\int_0^t||h'(s)||\d s )^2\d t,$$ and
\beg{ews}
&2L\int_0^T||y_t||_\infty||h'||\d t = 2L\int_0^Te^{t(L+a)}\int_0^t||h'(r)||\d r||h'(t)||\d t\\
&\leq LTe^{L(T+a)}\int_0^T||h'(r)||^2\d r-L(L+a)\int_0^Te^{t(L+a)}(\int_0^t||h'(s)||\d s )^2\d t
\end{ews}
Considering we have assume $a\geq0$, therefore
$$\int_0^T||y'(t)||^2\d t  \leq e^{2Ta^+ + T^2(a^+)^2e^{2Ta^+}}\Big(1+LTe^{T(L+a^+)}\Big)\int_0^T||h'(r)||^2\d r.$$
\end{proof}

\beg{prop}\label{prop&logsobolev}
Assume that $\sigma\equiv I$, $||\nabla F||\leq L$ and $\nabla F$ is uniformly continuous on bounded set of $\mathscr{C}\times\mathscr{C}$, then we have the following log Sobolev inequality holds
$$\E G^2\log G^2-\E G^2\log\E G^2\leq  2e^{2Ta^+ + T^2(a^+)^2e^{2Ta^+}}\Big(1+LTe^{T(L+a^+)}\Big) \mathscr{E}^\xi(G),$$
for all $G\in \mathscr{F}C_b^\infty(W^\xi(T))$. In particularly, for $G(\gamma)=g(\gamma(T))$, $g\in C^1_b(H)$, we have
$$\E (g^2\log g^2)-\E g^2\log\E g^2\leq 2e^{2Ta^+ + T^2(a^+)^2e^{2Ta^+}}\Big(1+LTe^{T(L+a^+)}\Big)\E ||\nabla g||^2(x^\xi(T)).$$
\end{prop}
\beg{proof}
Consider the gradient $(Dx)^*\nabla G$. Then for all adapted $h\in L^2(\Omega\ra W^{1,2}([0,T],H);\PP)$ with $||h||_{W^{1,2}}$ bounded $\PP$-a.s., we have
\beg{ews}\label{equ10}
&\E\<(Dx)^*\nabla G,h\>=\E\<\nabla G,D_h x\>=\E D_h G(x^\xi([0,T]))=\E G(x^\xi([0,T]))\int_0^T\<h'(t),\d W(t)\>.
\end{ews}
By martingale representation theorem, It's standard that  $$\E\Big[G(x^\xi([0,T]))|\mathscr{F}_t\Big]=\E G(x^\xi([0,T])) + \int_0^t\<\E\Big[(\ff \d {\d s} (Dx)^*\nabla G)|\mathscr{F}_s\Big],\d W(s)\>.$$
Let $m_t=\E\Big[G^2(x^\xi([0,T]))|\mathscr{F}_t\Big]$. By It'\^{o}'s formula, we have
\beg{ews}
&\E m_T\log m_T-m_0\log m_0=\int_0^T\E\ff {\Big(\E\Big[\ff {\d} {\d t}(Dx)^*\nabla G^2|\mathscr{F}_t\Big]\Big)^2}{2m_t}\d t\\
&\leq 2\int_0^T\E\ff{\Big(\E\Big[G\ff {\d} {\d t} (Dx)^*\nabla G|\mathscr{F}_t\Big]\Big)^2}{m_t}\d t\leq 2\E||(Dx)^*\nabla G||^2\leq  2 (1+LTe^{LT}) \E||\nabla G||^2.
\end{ews}
\end{proof}

\bigskip

Next, we shall extend the above result to the path space of segment processes. Let $W^\xi(T)$ the path space of $x^\xi(\cdot)$ on $[0,T]$, $\mathscr{H}$ consist of all the segment functions of $W^{1,2}_0([0,T],H)$, i.e.  $$\mathscr{H}=\{\psi_\cdot\ |\ \psi(\cdot)\in W^{1,2}_0([0,T],H)\ \mbox{extended to}\ [-\tau,0]\ \mbox{by zero}\},$$ $S$ be the natural embedding from $W^\xi(T)$ to $W^\xi_T$, i.e. $(S\gamma)_\cdot=\gamma_\cdot$ . Then we can introduce a inner product structure on $\mathscr{H}$ as follow such that it becomes a Hilbert space
$$\<\phi,\psi\>_\mathscr{H}:= \<S^{-1}\phi,S^{-1}\psi\>_{W^{1,2}_0},\ \phi,\psi\in \mathscr{H}.$$ Let $G\in C_b^1(W^\xi_T)$ and $\nabla_\eta G(\gamma)$ be the derivative of $G$ at $\gamma\in W^\xi_T$ along the direction $\eta\in \mathscr{H}$. If $\nabla_\cdot G(\gamma)$ gives a continuous linear functional on $\mathscr{H}$, we define the gradient $\nabla G(\gamma)$ as its Riesz representation, i.e. $\<\nabla G(\gamma),\eta\>_\mathscr{H}=\nabla_\eta G(\gamma)$. Let $W^{1,2}_\tau=W^{1,2}([-\tau,0],H)$. On $W^{1,2}([-\tau,0],H)$ we rig the inner product $\<\phi,\psi\>_{W^{1,2}_\tau}:=\<\phi(0),\psi(0)\>+\int_{-\tau}^0\<\dot{\phi}(s),\dot{\psi}(s)\>\d s,\ \phi,\psi\in W^{1,2}([-\tau,0],H)$.  Let $G(\gamma)=g(\gamma_{t_1},\cdots,\gamma_{t_n}),\ g\in C^1_b(\mathscr{C}^n)$, $\nabla^i g$ be the partial derivative of the i-th component, as an element of $W_\tau^{1,2}$ just as above. Then for any $\eta\in \mathscr{H}$
\beg{ews}
&\nabla_\eta G(\gamma)=\sum_{i=1}^{n}\<\nabla^i g(\gamma),\eta_{t_i}\>_{W^{1,2}_\tau}\\
&=\sum_{i=1}^{n}\Big[\int_{-\tau}^0\<\nabla^i g(\gamma)(s), \dot{\eta}(t_i+s)\>\d s+\<\nabla^i g(\gamma)(0),\eta(t_i)\>\Big]\\
&=\sum_{i=1}^{n}\Big[\int_0^T\<1_{[t_t-\tau,t_i]}(s)\nabla^i g(\gamma)(s-t_i),\dot{\eta}(s)\>\d s+\int_0^T\<1_{[s<t_i]}\nabla^i g(\gamma)(0),\dot{\eta}(s)\>\d s \Big]\\
&=\int_0^T\<\sum_{i=1}^{n}\Big[1_{[t_t-\tau,t_i]}(s)\nabla^i g(\gamma)(s-t_i)+1_{[s<t_i]}\nabla^i g(\gamma)(0)\Big],\dot{\eta}(s)\>\d s,
\end{ews}
we have
$$\Big(\ff {\d } {\d t}S^{-1}\nabla G(\gamma)\Big)(s)=\sum_{i=1}^{n}\Big[1_{[t_t-\tau,t_i]}(s)\nabla^i g(\gamma)(s-t_i)+1_{[s<t_i]}\nabla^i g(\gamma)(0)\Big].$$
A counterpart of $(Dx)^*\nabla G(\gamma)$ as in Proposition \ref{prop&logsobolev} is $S(Dx)^*S^{-1}\nabla G(\gamma)$. By these definition, just as in the previous discussion, we have the results on $W^\xi_T$.

\beg{prop}
Under the same assumption of Proposition \ref{prop&logsobolev}. For all $$\eta\in W^{1,2}_0([0,T],H)\bigcap L^2([0,T],\Dom(A)),$$ we have the integration by part formula
$$\E\nabla_{S\eta}G(\gamma)=\int_0^T G(x^\xi_{[0,T]})\<\dot{\eta}(t)+A\eta(t)-\nabla_{\eta_t}F(x_t),\d W(t)\>.$$ Let $\mathscr{E}_T^\xi(\Phi,\Psi)=\E\<\nabla\Phi,\nabla\Psi\>_{\mathscr{H}}(x^\xi_{[0,T]})$. Then $\Big(\mathscr{E}_T^\xi,\mathscr{F}C_b^\infty(W^\xi_{T})\Big)$ is closable in $L^2(W^\xi_{T},\Pi^\xi_{T})$. Log Sobolev inequality holds
$$\E G^2\log G^2-\E G^2\log \E G^2\leq 2e^{2Ta^+ + T^2(a^+)^2e^{2Ta^+}}\Big(1+LTe^{T(L+a^+)}\Big)\E ||\nabla G||_\mathscr{H}^2.$$
In particularly, for $G(\gamma)=g(\gamma_T)$, $g\in C_b^1(\mathscr{C})$,
$$P_Tg^2\log g^2(\xi)-P_Tg^2(\xi)\log P_Tg(\xi)\leq 2(T+1)e^{2Ta^+ + T^2(a^+)^2e^{2Ta^+}}\Big(1+LTe^{T(L+a^+)}\Big)\E ||\nabla g||^2_{W_\tau^{1,2}}(\xi).$$
\end{prop}

\section{Stochastic evolution equation with non-Lipschitz coefficients}
Here we consider the following equation in Hilbert space $H$
\be\label{equ1}
\d x(t)=-Ax(t)\d t +B(x(t))\d t+Q\d W(t).
\ee
We shall use the notation following
\begin{eqnarray}
(V_\theta,||\cdot||_{V_\theta})=(\Dom(A^\frac{\theta}{2}),||A^\frac{\theta}{2}\cdot||),\ ||\cdot||_Q=||Q^{-1}\cdot||.
\end{eqnarray}
The coefficients of the equation may satisfy some of the following conditions
\begin{enumerate}[({A}1)]
\item $A$ is positive self adjoint operator with $A\geq\lambda_0>0$, $Q\in\mathscr{L}_{HS}(H)$ is non-degenerated,\label{A1}
\item $B$ is hemicontinuous, i.e. the map $s\ra \<B(v_1+sv_2),v\>$ is continuous on $\mathbb{R}$, and there exists $\gamma\in[0,2)$, $\alpha\in[0,1]$ and $K_1,K_2\geq 0$ such that
\begin{eqnarray}
\<B(u)-B(v),u-v\>&\leq&(\rho(v)+K_1)||u-v||_V^\gamma||u-v||^{2-\gamma}\label{equ7}\\
\<B(u-v),v\>&\leq& K_2||v||_{V_\alpha}||u-v||_V^{\gamma}||u-v||^{2-\gamma}\label{equ8}
\end{eqnarray}
where $\rho:V\ra\mathbb{R}^+$ is measurable, locally bounded function, $\rho(0)=0$,\label{A2}
\item There exists $\theta\in(0,1]$ and $K_3>0$ such that \label{A3}
    \be
    ||u||_Q^2\leq K_3||u||^2_{V_{\theta}},
    \ee
\item There exists a constant $K_4>0$ such that
\be
||B(u)-B(v)||^2_Q\leq \beta(u-v)(1+||v||_V+||u||_V)^2,
\ee
where $\beta: V\ra \mathbb{R}^+$ is locally bounded measurable function. \label{A4}
\item $B$ is G\^{a}teaux differentiable from $V$ to $Q(H)$ and there exists $K_4\geq 0$  such that \label{A5}
\be
||\nabla B(v)||_Q\leq K_4(1+||v||_V),
\ee
here we endow $Q(H)$ with the norm $||\cdot||_Q$ such that it becomes a Banach space,
\item There is $K_5\geq 0$ such that \label{A6}
\be
\<B(w),w\>\leq K_5(1+||w||^2), \ \forall w\in V.
\ee
\end{enumerate}

\begin{remark}\label{rm2}
(1) By (\ref{equ8}), we have
\be
\<B(w),v\>=\<B(w+v-v),v\>\leq K_2||v||_{V_\alpha}||w||_V^\gamma||w||^{2-\gamma},\ \forall v,w\in V,
\ee
thus
\be
||B(w)||_{V^*}\leq C||w||_V^\gamma||w||^{2-\gamma},\ \forall w\in V,
\ee
and from (\ref{equ7}),
\be
\<B(w),w\>\leq ||B(0)||\cdot||w||+K_1||w||_V^\gamma||w||^{2-\gamma}, \ \forall w\in V.
\ee
Therefore, by \cite{RoLiu2010} and directly calculus  we can prove that under the conditions (A\ref{A1}) and (A\ref{A2}), equation (\ref{equ1}) has uniqueness strong solution.

(2) It's easy to see that (A\ref{A5}) implies that (A\ref{A4}) holds in the following form
\be
||B(u)-B(v)||_Q\leq 2K_4||u-v||_V(1+||u||_V+||v||_V).
\ee

(3) Though it's easy to see that Navier-Stokes operator satisfies (A\ref{A2}), but unfortunately, it does not satisfies (A\ref{A3}) to (A\ref{A5}).
\end{remark}

\begin{theorem}\label{theorem2}
Assume that (A\ref{A1}) to (A\ref{A4}) hold and $e\in\Dom(A^{\frac{1+\theta}{2}})$, then the shift log-Harnack inequality holds
\be
P_T\log f(x)\leq \log P_Tf(e+\cdot)(x) + \Psi_2(T,e),\ \ \forall f\in \mathscr{B}_b(H),
\ee
here
\begin{ews}
\Psi(x,T,e)=&\ C\Big\{\ff {(T+1)} T ||A^{(1+\theta)/2}e||^2 +b_e + b_e \Big((||Q||^2_{HS}+||B(0)||^2+||x||)T\\&+\ff {||e||^2} 2+ \ff {\sqrt{2-\gamma}} 4 ||A^{1/4}e||^2||A^{(2\alpha+\gamma-2)/4}e||^{2/(2-\gamma)}\Big)\\
&\quad\times\exp\Big\{CT\Big[1+\ff {2-\gamma} 4 ||A^{(2\alpha+\gamma-2)/4}e||^{4/(2-\gamma)}\Big]\Big\} \Big\},
\end{ews}
where C is constant depending on $\gamma,\ K_1,\ K_2, K_3$ and
\be
b_e:=\sup_{||v||_V\leq ||e||_V}{\beta(v)}.
\ee
If we assume that (A\ref{A1}) to (A\ref{A3}) and (A\ref{A5}) hold, then for all $\phi\in L^2([0,T],V_\theta),\ L^I_T\phi=e$, the integration by part formula holds
\be
P_T\nabla_e f(x)=\E f(x(T))\int_0^T\<Q^{-1}(\phi(t)-\nabla_{\Gamma(t)}B(x(t))),\d W(t)\>, \ \forall f\in C_b^1(H),
\ee
where
\be
\Gamma(t)=\int_0^te^{-(t-s)A}\phi(s)\d s.
\ee
\end{theorem}
\noindent\emph{Proof.}
Consider the operator $L^{A^{\theta/2}}_T$ which maps from $L^2([0,T],V_\theta)$ to $H$. Since $A^{-\frac{\theta}{2}}: H\ra V_\theta$ is isometric,
\be
L^{A^{\theta/2}}_T: L^2([0,T],V_\theta)\ra V_1
\ee
is surjective, by proposition \ref{prop1}. Note that $A$ is self adjoint, thus $A^{-\frac{\theta}{2}}$ is adjoint of $A^{\frac{\theta}{2}}$ as an operator from $V_\theta$ to $H$, then $\mbfR^{A^{\theta/2}}_T=\int_0^Te^{-2tA}\d t$. Firstly we shall  choose special $\phi$ to get log Haranck inequality. Since $A^{\theta/2}e\in V_1$, as in Lemma \ref{lmm3}, replacing $x$ by $A^{\theta/2}e$ and $B^{-1}$ by $A^{-\theta/2}e$, we have $\phi(t)=\ff 1 T e^{-(T-t)A}e + \ff {2t} T Ae^{-(T-t)A}e$. Then
\begin{eqnarray}
\Gamma(t)&=& \ff t T e^{-(T-t)A}e,\ L_T\phi=e,\\
\int_0^T||A^{\frac{\theta}{2}}\phi(s)||^2\d s&=&\int_0^T||\phi(s)||_{V_\theta}^2\d s\leq \ff {2(1+T)} T||A^{\ff {1+\theta} 2}||^2,\label{inequ4}
\end{eqnarray}
and by (A\ref{A3})
\be
\int_0^T||Q^{-1}\phi(s)||^2\d s\leq K_3\int_0^T||A^{\frac{\theta}{2}}\phi(s)||^2\d s\leq \ff {2(1+T)K_3} T||A^{\ff {1+\theta} 2}||^2.
\ee
We construct another process
\be
\d y(t)=-Ay(t)\d t+B(x(t))\d t+Q\d W(t)+\phi(t)\d t,\ y(0)=x,
\ee
then $y(t)=x(t)+\Gamma(t)$, in particular, $y(T)=x(T)+e$. Let
\be
\d \tilde{W}(t)=\d W(t)+Q^{-1}\phi(t)\d t+Q^{-1}(B(x(t))-B(x(t)+\Gamma(t)))\d t,
\ee
and
\begin{ews}
R_t=&\exp\left[-\int_0^t\<Q^{-1}(\phi(s)+B(x(s))-B(x(s)+\Gamma(s))),\d W(s)\>\right.\\
&\left.-\frac{1}{2}\int_0^t||Q^{-1}(\phi(s)+B(x(s))-B(x(s)+\Gamma(s)))||^2\d s\right],
\end{ews}
we can rewrite $y$ as
\be
\d y(t)=-Ay(t)\d t +B(y(t))\d t +Q\d \tilde{W}(t).
\ee
Next we shall prove that $\{\tilde{W}(t)\}_{t\in[0,T]}$ is $R_T\mathbb{P}$-Brownian Motion, then $y$ is a weak solution of equation (\ref{equ1}), and since equation (\ref{equ1}) has pathwise unique solution, $y$ and $x$ has the same law under the probability measures respectively, then
\be
\E f(x(T))=P_Tf(x)=\E R_Tf(y(T))=\E R_Tf(x(T)+e),
\ee
therefore the argument in \cite{Wang2012} can be applied.
To this end, we shall adapt the argument in \cite{Wang2011,Wang2012} to estimate $\E R_t\log{R_t}$.
Note that
\begin{eqnarray}
\sup_{t\in[0,T]}||\Gamma(t)||_V\leq ||e||_V,\ \sup_{t\in[0,T]}\beta(\Gamma(t))\leq \sup_{||v||_V\leq ||e||_V}{\beta(v)}.\label{inequ1}
\end{eqnarray}
Then 
\begin{ews}
||Q^{-1}(B(x(t))-B(x(t)+\Gamma(t)))||^2&\leq \beta(\Gamma(t))(1+||x(t)||_V+||x(t)+\Gamma(t)||_V)^2\\
&\leq 3b_e(1+3||x(t)||^2_V+2||\Gamma(t)||^2_V)
\end{ews}
Let
\be
\tau_n=\inf\{t\in[0,T]\ |\int_0^t||x(s)||_V^2\d s+||x(t)||^2\geq n \}.
\ee Then by Girsanov theorem, for $s\leq T$, $\{\tilde{W}(t)\}_{t\leq s\wedge\tau_n}$ is Brownian Motion under the probability $R_{s\wedge\tau_n}\mathbb{P}$. Rewrite the equation of $x$, we have
\begin{ews}\label{equ2}
\d x(t)&=-Ax(t)\d t +B(x(t))\d t+Q\d W(t)\\
&=-Ax(t)\d t+B(x(t))\d t+Q\d \tilde{W}(t)-\phi(t)\d t-(B(x(t))-B(x(t)+\Gamma(t)))\d t\\
&=-Ax(t)\d t+B(x(t)+\Gamma(t))\d t+Q\d \tilde{W}(t)-\phi(t)\d t,\ t\leq s\wedge\tau_n,
\end{ews}
by It'\^{o}'s formula and (A\ref{A2}), as what we do to equation (\ref{equ2}), we rewrite it in the form of $\tilde{W}$, then get that, for any $t\leq s\wedge \tau_n$
\begin{ews}
&\d ||x(t)||^2+2||x(t)||^2_V\d t-||Q||_{HS}\d t+2\<\phi(t),x(t)\>\d t\\
=&\ 2\<B(x(t)+\Gamma(t)),x(t)\>\d t+2\<Q\d \tilde{W}(t),x(t)\>\\
=&\ 2\<B(x(t)+\Gamma(t)),x(t)+\Gamma(t)\>\d t-2\<B(x(t)+\Gamma(t)),\Gamma(t)\>\d t+2\<Q\d \tilde{W}(t),x(t)\>\\
\leq &\ \lb 2||B(0)||\cdot||x(t)+\Gamma(t)||+2K_1||x(t)+\Gamma(t)||_V^\gamma||x(t)+\Gamma(t)||^{2-\gamma}\rb\d t\\
&+2K_2||\Gamma(t)||_{V_\alpha}||x(t)+\Gamma(t)||^\gamma_V||x(t)+\Gamma(t)||^{2-\gamma}\d t+2\<Q\d \tilde{W}(t),x(t)\>.
\end{ews}
In following $C$ is constant depend on $\gamma,\ K_1,\ K_2$ may change from line to line. By B-D-G inequality and H\"{o}lder inequality, we have
\begin{ews}
&\E R_{s\wedge\tau_n}\sup_{r\in[0,t\wedge\tau_n]}||x(r)||^2+\E R_{s\wedge\tau_n}\int_0^{t\wedge\tau_n}||x(r)||^2_V\d r\\
\leq &(C+||B(0)||^2+||x||+||Q||_{HS}^2)t+2\E R_{s\wedge\tau_n}\sup_{r\in[0,t\wedge\tau_n]}\left|\int_0^r\<Q\d \tilde{W}(u),x(u)\>\right|\\
&+\int_0^t||\phi(r)||^2\d r+C\int_0^T||\Gamma(t)||_V^2(1+||\Gamma(t)||_{V_\alpha}^\frac{2}{2-\gamma})\d t\\
&+C\E R_{s\wedge\tau_n}\lb \int_0^{t\wedge\tau_n}||x(r)||^2\d r +\int_0^{t\wedge\tau_n}||x(r)||^2||\Gamma(r)||^\frac{2}{2-\gamma}_{V_\alpha}\d r\rb,\\
\leq &(C+||Q||^2_{HS}+||B(0)||^2+||x||)t+\int_0^t||\phi(r)||^2\d r\\
&+C\int_0^T||\Gamma(t)||_V^2(1+||\Gamma(t)||_{V_\alpha}^\frac{2}{2-\gamma})\d t\\
&+C\E R_{s\wedge\tau_n}\lb \int_0^{t\wedge\tau_n}||x(r)||^2\d r +\int_0^{t\wedge\tau_n}||x(r)||^2||\Gamma(r)||^\frac{2}{2-\gamma}_{V_\alpha}\d r\rb,
\end{ews}
In order to use the Gronwall's lemma, we need more calculate. Note that for the last term, we have
\begin{ews}
&\E R_{s\wedge\tau_n}\int_0^{t\wedge\tau_n}||x(r)||^2||\Gamma(r)||^\frac{2}{2-\gamma}_{V_\alpha}\d r\\
\leq &C\lb\int_0^t||\Gamma(r)||_{V_\alpha}^{\frac{4}{2-\gamma}}\d r\rb^\frac{1}{2}\E R_{s\wedge\tau_n}\lb\int_0^{t\wedge\tau_n}||x(r)||^4\d r\rb^\frac{1}{2}\\
\leq & C\lb\int_0^t||\Gamma(r)||_{V_\alpha}^{\frac{4}{2-\gamma}}\d r\rb^\frac{1}{2}\E R_{s\wedge\tau_n}\lb\sup_{r\in[0,t\wedge\tau_n]}{||x(r)||}\rb\lb\int_0^{t\wedge\tau_n}||x(r)||^2\d r\rb^\frac{1}{2}\\
\leq & \frac{1}{2}\E R_{s\wedge\tau_n}\sup_{r\in[0,t\wedge\tau_n]}{||x(r)||^2} +C\lb\int_0^t||\Gamma(r)||_{V_\alpha}^{\frac{4}{2-\gamma}}\d r\rb\E R_{s\wedge\tau_n}\int_0^{t\wedge\tau_n}||x(r)||^2\d r\\
\end{ews}
In order clear relation with $e$, we shall calculate the integration term relate to $\Gamma$. By Minkowski inequality
\begin{ews}\label{inequ3}
\int_0^T||\Gamma(r)||_{V_\alpha}^{\frac{4}{2-\gamma}}\d r
&=T^{-4/(2-\gamma)}\int_0^T||rA^{\alpha/2}e^{-(T-r)A}e||^{4/(2-\gamma)}\d r\\
&=T^{-4/(2-\gamma)}\Big(\int_0^T\Big(\int_{\lambda_0}^\infty r^2\lambda^\alpha e^{-2(T-r)\lambda}\d ||E_\lambda e||^2\Big)^{2/(2-\lambda)}\Big)^{\ff 2 {2-\gamma} \ff {2-\gamma} 2}\\
&\leq T^{-4/(2-\gamma)}\Big(\int_{\lambda_0}^\infty \lambda^\alpha\Big(\int_0^T r^{4/(2-\gamma)}e^{-4(T-r)\lambda/(2-\lambda)}\d r\Big)^{(2-\gamma)/2}\d ||E_\lambda e||^2\Big)^{2/(2-\gamma)}\\
&\leq \ff {2-\gamma} {4} ||A^{(2\alpha+\gamma-2)/4}e||^{4/(2-\gamma)}.
\end{ews}
In particular, for $\alpha=1$ and $\gamma=0$, we have
\be\label{inequ5}
\int_0^T||\Gamma(r)||_{V}^{2}\d r\leq \ff {||e||^2} 2.
\ee
For $\alpha=1,\ \gamma=1$,
\beg{ews}
\int_0^T||\Gamma(r)||_{V}^{2}||\Gamma(r)||^\frac{2}{2-\gamma}_{V_\alpha}\d r&\leq\lb\int_0^T||\Gamma(r)||_{V}^{4}\d r\rb^\frac{1}{2}\lb\int_0^T||\Gamma(r)||^\frac{4}{2-\gamma}_{V_\alpha}\d r\rb^\frac{1}{2}\\
&\leq \ff {\sqrt{2-\gamma}} 4 ||A^{1/4}e||^2||A^{(2\alpha+\gamma-2)/4}e||^{2/(2-\gamma)}.
\end{ews}
At last  
\be
\int_0^T||\phi(r)||^2\d r\leq \ff {2(T+1)} {T}||A^{1/2}e||^2.
\ee
Therefore
\begin{ews}
&\E R_{s\wedge\tau_n}\sup_{r\in[0,t\wedge\tau_n]}||x(r)||^2+\E R_{s\wedge\tau_n}\int_0^{t\wedge\tau_n}||x(r)||^2_V\d r\\
\leq &C\Big\{ (||Q||^2_{HS}+||B(0)||^2+||x||)T+\ff {||e||^2} 2+ \ff {\sqrt{2-\gamma}} 4 ||A^{1/4}e||^2||A^{(2\alpha+\gamma-2)/4}e||^{2/(2-\gamma)}\Big\}\\
&+C\lb 1+\ff {2-\gamma} 4 ||A^{(2\alpha+\gamma-2)/4}e||^{4/(2-\gamma)}\rb\int_0^{t\wedge\tau_n}||x(r)||^2\d r
\end{ews}
By Gronwall's inequality
\begin{ews}
\E \sup_n\Big[R_{s\wedge\tau_n}\sup_{r\in[0,t\wedge\tau_n]}||x(r)||^2+\E R_{s\wedge\tau_n}\int_0^{t\wedge\tau_n}||x(r)||^2_V\d r\Big] \leq \infty.
\end{ews}
By these estimate and (A\ref{A4}), we have
\begin{ews}
&\E R_{s\wedge\tau_n}\log{R_{s\wedge\tau_n}}
= \ \frac{1}{2}\E R_{s\wedge\tau_n}\int_0^{s\wedge\tau_n}||Q^{-1}(\phi(t)+B(x(t))-B(x(t)-\Gamma(t)))||^2\d t\\
\leq &\ K_3\int_0^T||\phi(t)||^2_{V_\theta}\d t + 3b_e\E R_{s\wedge\tau_n}\int_0^{s\wedge\tau_n}\lb1+4||x(t)||^2_V+||\Gamma(t)||^2_V\rb\d t \\
=&\ \Psi(x,T,e).
\end{ews}
therefore, as in \cite{Wang2012}, we can prove that $\{\tilde{W}(t)\}_{t\in[0,T]}$ is B.M. and
\be
\E R_T\log{R_T}\leq \Psi(x,T,e).
\ee
By this estimate and Young's inequality, we have the shift log-Harnack inequality,
\begin{ews}
&P_T\log f(x)= \E_\mathbb{Q}\log f(y_T)=\E R_T\log f(x_T+e)\\
\leq &\E R_T\log R_T+\log \E f(x_T+e)\leq \log P_T f(e+\cdot)(x)+\Psi(x,T,e).
\end{ews}
For Integration by part formula, one can choose any $\phi\in L^2([0,T],V_\theta)$ such that $L^I_T\phi=e$. Replacing $e$ by $\epsilon e$, $\phi $ by $\epsilon \phi$ and $\Gamma$ by $\epsilon \Gamma$, just as in the case $\epsilon=1$ above, and by Lemma \ref{lmm2}, we have
\be
\frac{\d }{\d \epsilon}|_{\epsilon=0}R^\epsilon_T=-\int_0^T\<Q^{-1}\Big(\phi(t)-\nabla_{\Gamma(t)}B(x(t))\Big),\d W(t)\>,
\ee
holds in $L^1(\mathbb{P})$, then
\be
P_T\nabla_e f(x)=\E f(x(T))\int_0^T\<Q^{-1}\Big(\phi(t)-\nabla_{\Gamma(t)}B(x(t)))\Big),\d W(t)\>, \ f\in C_b^1(H).
\ee

\qed

We adapted the argument in \cite{GuiWang2012} to prove that
\begin{lemma}\label{lmm2}
Under conditions (A\ref{A1}) to (A\ref{A3}) and (A\ref{A5}), then
$\{\frac{|R^\epsilon-1|}{\epsilon}\}_{\epsilon\in(0,1)}$
is uniformly integrable w.r.t $\mathbb{P}$, consequently
\be
\frac{\d }{\d \epsilon}|_{\epsilon=0}R^\epsilon_T=-\int_0^T\<Q^{-1}\Big(\phi(t)-\nabla_{\Gamma(t)}B(x(t))\Big),\d W(t)\>,
\ee
holds in $L^1(\mathbb{P})$.
\end{lemma}
\noindent\emph{Proof.}
Denote
\begin{eqnarray}
\Theta_1^\epsilon(s)&=&Q^{-1}(\epsilon\phi(s)+B(x(s))-B(x(s)+\epsilon\Gamma(s)))\\
\Theta_2^\epsilon(s)&=&Q^{-1}(\phi(s)+\nabla_{\Gamma(s)}B(x(s)+\epsilon\Gamma(s)).
\end{eqnarray}
Since
\begin{ews}
&\E \int_0^{T}\sup_{\epsilon\in[0,1)}||\nabla_{\Gamma(r)}B(x(r)+\epsilon\Gamma(r))||_Q^2\d r\\
\leq &\E \int_0^{T}||\Gamma(s)||^2_V\lb||x(s)||_V^2+||\Gamma(s)||_V^2\rb\d s\\
\leq &\sup_{s\in[0,T]}||\Gamma(s)||_V^2\E \int_0^{T}||x(s)||_V^2\d s +\int_0^T||\Gamma(s)||^4_V\d s<\infty,
\end{ews}
we have, for any $\epsilon\in[0,1)$,
\begin{ews}
\frac{\d }{\d \epsilon}R^\epsilon_{T}=&-R^\epsilon_{T}\int_0^{T}
\<\Theta_2^\epsilon(s),\d W(s)\>-R_{T}^\epsilon\int_0^{T}\<\Theta_1^\epsilon(s),\Theta_2^\epsilon(s)\>\d s,\ a.s.,
\end{ews}
and then
\begin{ews}
\frac{|R_{T}^\epsilon-1|}{\epsilon}=&\left|\frac{1}{\epsilon}\int_0^\epsilon R_{T}^r\lb\int_0^{T}
\<\Theta_2^r(s),\d W(s)\>+\<\Theta_1^r(s),\Theta_2^r(s)\>\d s\rb\d r\right|\\
\leq & \left|\frac{1}{\epsilon}\int_0^\epsilon R_{T}^r\int_0^{T}\<\Theta_2^r(s),\d W(s)\>\d r\right|+\left|\frac{1}{\epsilon}\int_0^\epsilon R^r_{T}\int_0^{T}\<\Theta_1^r(s),\Theta_2^r(s)\>\d s\d r\right|.
\end{ews}
Note that
\begin{ews}
&\left|\frac{1}{\epsilon}\int_0^\epsilon R^r_{T}\int_0^{T}\<\Theta_1^r(s),\Theta_2^r(s)\>\d s\d r\right|\\
\leq & \frac{1}{\epsilon}\int_0^\epsilon R_{T}^r\int_0^{T}r\left[ ||\phi(s)||_Q+||\Gamma(s)||_V(||\Gamma(s)||_V+||x(s)||_V)\right]^2\d s\d r\\
\leq &\int_0^1 R_{T}^r\int_0^{T}\left[ ||\phi(s)||_Q+||\Gamma(s)||_V(||\Gamma(s)||_V+||x(s)||_V)\right]^2\d s\d r
\end{ews}
and just as in the case of $\epsilon=1$, we can prove that
\begin{ews}
&\E \int_0^1 R_{T}^r\int_0^{T}\left[ ||\phi(s)||_Q+||\Gamma(s)||_V(||\Gamma(s)||_V+||x(s)||_V)\right]^2\d s\d r\\
\leq &\int_0^{T} ||\phi(s)||^2_Q\d s+\int_0^T||\Gamma(s)||^4_V\d s +\sup_{s\in[0,T]}||\Gamma(s)||_V^2\int_0^1\E R_{T}^r\int_0^T||x(s)||_V^2\d s\d r\\
\end{ews}
By these estimate, follow the line of \cite[Lemma 2.4.]{GuiWang2012} completely, one can prove the lemma.
\qed

\begin{corollary}
Assume that (A\ref{A1}) to (A\ref{A3}) and (A\ref{A6}) hold, further more (A\ref{A5}) or (A\ref{A4}) hold in the following form
\be\label{inequ6}
||B(u)-B(v)||_Q\leq K_4||u-v||(1+||u||_V+||v||_V),
\ee
let
\be
\delta_e=\frac{e^{(\lambda_0-2K_5)^-T}}{18K_4||Q||^2||e||^2_VT}
\ee
then for $r\in(0,\sqrt{\delta_e})$ and $p>\frac{\sqrt{8\delta_e r^2+r^4}+2\delta_e+r^2}{2\delta_e-r^2}$,
\begin{ews}
(P_Tf(x))^p\leq &\ P_Tf^p(re+\cdot)(x)\exp\Big\{\frac{p-1}{4||Q||^2}\Big[\frac{||x||^2}{T}+\Big(||Q||^2_{HS}+2K_5\Big)[(\lambda_0-2K_5)^+T\vee1]\Big]\\
&\ +\ff {(p+1)p} {2(p-1)} \Big[\ff {2K_3(T+1)} {T} ||A^{(1+\theta)/2}||^2 + \ff 3 2 ||A^{1/4}e||^4+\ff {3K_4} {2} ||e||^2\Big].
\end{ews}
If strengthen (A\ref{A4}) to be
\be
||B(u)-B(v)||_Q\leq \beta(u-v),
\ee
then, for any $p>1$, the following shift Harnack inequality holds
\be
(P_Tf)^p\leq P_Tf^p(e+\cdot)\exp\Big[\frac{p}{p-1}\Big(\frac{2||A^\frac{1+\theta}{2}e||^2}{1-e^{-2\lambda_0T}}+\int_0^T\beta(\Gamma(s))\d s\Big)\Big]
\ee
\end{corollary}
\noindent\emph{Proof.}
By Remark \ref{rm2}, we assume (\ref{inequ6}) holds. We adapted the technology used in \cite[Lemma 3.1]{WangXu2012}. Let $\lambda=\frac{c_0}{4||Q||^2}$ and $$\beta(t)=\Big[\frac{c_0(1-e^{-(\lambda_0-2K_5)t})}{\lambda_0-2K_5}+c_0Te^{-(\lambda_0-2K_5)t}\Big]^{-1},$$
for $\lambda_0-2K_5=0$ we define it as $\frac{1}{c_0T}$. By It'\^{o}'s formula for $||x(t)||^2\beta(t)$ and H\"{o}lder inequality, we can prove that
\begin{ews}
&\E \exp\Big[2\lambda\int_0^{T\wedge\tau_n}\beta(t)||x(t)||_V^2\d t-\lambda\beta(0)||x||^2-\lambda{||Q||_{HS}^2+2K_5}\int_0^T\beta(t)\d t\Big]\\
&\leq \Big[\E\exp\Big[2\lambda\int_0^{T\wedge\tau_n}[(2K_5-\lambda_0)\beta(t)+\beta'(t)]\cdot||x(t)||^2\d t+4\lambda\int_0^{T\wedge\tau_n}\beta(t)||x(t)||^2\d t\Big]\Big]^{1/2}\\
&\quad\times\Big[\E\exp[2\lambda\int_0^{T\wedge\tau_n}\beta(t)||x(t)||^2_V\d t]\Big]^{1/2},
\end{ews}
then by the definition of $\lambda$ and $\beta$, we have
\beg{ews}
&\E\exp\Big[\frac{e^{(\lambda_0-2K_5)^-T}}{2T||Q||^2}\int_0^T||x(t)||_V^2\d t\Big]\\
&\leq \exp\Big\{\frac{||x||^2}{2T||Q||^2}+\frac{||Q||^2_{HS}+2K_5}{2||Q||^2}[(\lambda_0-2K_5)^+T\vee1]\Big\}.
\end{ews}
For $r\in(0,\sqrt\delta_e)$, just replacing $e$ by $re$ in Theorem \ref{theorem2}, for and $p>\frac{\sqrt{8\delta_e r^2+r^4}+2\delta_e+r^2}{2\delta_e-r^2}$, we can prove that
\begin{ews}
&(\E R_T^\frac{p}{p-1})^{p-1}\exp \left\{-\frac{p(p+1)}{2(p-1)}\lb K_3\int_0^T||\phi(t)||_{V_\theta}\d t+6K_4\int_0^T||\Gamma(t)||_V^4\d t + 3K_4\int_0^T||\Gamma(t)||_V^2\d t\rb\right\}\\
&\leq\lb\E\exp\left[\frac{18K_4p(p+1)}{(p-1)^2}\int_0^T||\Gamma(t)||_V^2||x(t)||_V^2\d t \right]\rb^\frac{p-1}{2}\\
&\leq\lb\E\exp\left[\ff {9K_4p(p+1)}{(p-1)^2}||re||_V^2
\int_0^T||x(t)||_V^2\d t \right]\rb^\frac{p-1}{2},
\end{ews}
by the definition of $\delta_e$, we have
\be
\frac{9K_4p(p+1)||e||_V^2r^2}{(p-1)^2}
\leq  \frac{e^{(\lambda_0-2K_5)^-T}}{2T||Q||^2},
\ee
then
\begin{ews}
&(\E R_T^\frac{p}{p-1})^{p-1}\exp \left\{-\frac{p(p+1)}{2(p-1)}\lb K_3\int_0^T||\phi(t)||_{V_\theta}^2\d t+6K_4\int_0^T||\Gamma(t)||_V^4\d t + 3K_4\int_0^T||\Gamma(t)||_V^2\d t\rb\right\}\\
&\leq \exp\Big\{\frac{p-1}{4||Q||^2}\Big[\frac{||x||^2}{T}+\Big(||Q||^2_{HS}+2K_5\Big)[(\lambda_0-2K_5)^+T\vee1]\Big]\Big\}.
\end{ews}
Combine this with (\ref{inequ3}), (\ref{inequ4}) and (\ref{inequ5}), we prove the first inequality. The second inequality is similar to corollary \ref{coro3}.

\qed

\begin{corollary}
For hyperdissipative stochastic Navier-Stokes/Burgers equation in \cite{WangXu2012}, (A\ref{A1}) to (A\ref{A6}) hold.
\end{corollary}
\noindent\emph{Proof.}
We just have to verify the conditions (A\ref{A1}) to (A\ref{A6}). (A\ref{A3}) is the (A0) there, and by the bilinear, it's Fr\'{e}chet differentiable form $V$ to $Q(H)$, and by (A3) in \cite{WangXu2012},
\be
||\nabla_uB(v)||_Q=||B(u,v)+B(v,u)||_Q\leq C||u||_{V_\theta}||v||_{V_\theta}\leq C||u||_V||v||_V,
\ee
then (A\ref{A5}) holds. By (A3) in \cite{WangXu2012},
\begin{ews}
||B(u)-B(v)||_Q=&||B(u-v+v)-B(v)||_Q\\
=&||B(u-v)+B(u-v,v)+B(v,u-v)||_Q \\
\leq &||B(u-v)||_Q+||B(u-v,v)||_Q+||B(v,u-v)||_Q\\
\leq &C\lb||u-v||^2_{V_\theta}+2||u-v||_{V_\theta}||v||_{V_\theta}\rb\\
\leq &C||u-v||_V\lb||u||_V+||v||_V\rb,
\end{ews}
then (A\ref{A4}) holds with $\beta(\cdot)=C||\cdot||_V$ and $K_4=0$. By (A2) in \cite{WangXu2012}, and
\begin{ews}
\<B(u)-B(v),u-v\>&=\<B(u-v+v)-B(v),u-v\>\\
&=\<B(u-v,v)+B(v,u-v)+B(u-v),u-v\>\\
&=\<B(v,u-v),u-v\>+\<B(u-v,v),u-v\>\\
&\leq C(||v||\cdot||u-v||_V||u-v||+||u-v||\cdot||v||_V||u-v||)\\
&\leq C||v||_V||u-v||_V||u-v||.
\end{ews}
(\ref{equ8}) and (\ref{equ7}) holds for $\gamma=1$, $\rho=||\cdot||_V$, the hemicontinuous follows from bilinear. Therefore, we prove the corollary.

\qed

At last, we give a simple corollary to discuss the density of solution of equation (\ref{equ1}).
\beg{corollary}
Under the condition of Theorem \ref{theorem2}. Let $\pi_n$ be the orthogonal projection from $H$ to some $n$ dimension subspace $H_n$. Then the distribution of $\pi_n x(T)$ has density $\rho_n$ with respective to Lebesgue measure on $H_n$ and $$\nabla\log\rho_n(x)=-\E(N|\pi_n x(T)=x),\ \ \PP_{\pi_n x(T)}\mbox{-a.s.},$$
where $N$ such that
$P_T\nabla_e f= \E f(x(T))N(e).$
\end{corollary}
\beg{proof}
For all $f\in C_b^1{H_n}$, let $f_n(x)=f(\pi_n x)\ \forall x\in H$. Then for all $h\in H_n$,
$$\E\nabla_h f(\pi_n x(T))=\E\nabla_h f_n(x(T))=\E f_n(x(T))N(h)=\E f(\pi_n x(T))N(h),$$
by Theorem 2.4 in \cite{Wang2012}, we prove the corollary.
\end{proof}
\beg{remark}
The more interesting case for the density of the projection of the solution is SPDE driven by degenerate noise, see \cite{BM2007} and reference there in.
\end{remark}

For more applications of shift Harnack inequality and integration by part formula, one can see \cite{Wang2012}.


\bigskip

\textbf{Acknowledgement}
Thanks Professor Feng-Yu Wang for useful suggestions. It's kind of Mr. Ouyang to email the author to mention their work \cite{OuyangRW2012} after the first version was submit.

\begin{thebibliography}{99}

\bibitem{BM2007}
    Bakhtin, Y. and Marttingly, J.~C. (2007).
    Malliavin Calculus for Infinite-dimensional systems with additive noise.
    \textit{J. Funct. Anal.}
    \textbf{249} 307--353.

\bibitem{BYWang}
    Bao, J., Yuan, C. and Wang, F.-Y. (2012).
    Bismut formulae and applications for functional SPDEs, arXiv:1110.5150.

\bibitem{BDPDM2007}
    A. Bensoussan, G. Da Prato, M. C. Delfour, S. K. Mitter. (2007).
    Representation And Control Of Infinite Dimensional Systems. 2nd ed.
    \textit{Birkh\"{a}user}.

\bibitem{BH2010}
    Boufoussi, B., Hajji, S. (2010).
    Successive approximation of neutral functional stochastic differential equations in Hilbert spaces.
    \textit{Ann. Math. Blaise Pascal}
    \textbf{17} 183--197.

\bibitem{DPZ1992}
    G. Da Prato, J. Zabczyk (1992).
    Stochastic Equations in Infinite Dimensions.
    \textit{Cambridge University Press.}

\bibitem{Dri1997}
    Driver, B. (1997).
    Integration by parts for heat kernel measures revisited.
    \textit{J. Math. Pures Appl.}
    \textbf{76}  703--737.

\bibitem{ERS2009}
    Es-Sarhir, A., Von Renesse, M. and Scheutzow, M. (2009).
    Harnack inequality for functional SDEs with Bounded Memory.
    \textit{Elect. Comm. in Probab.}
    \textbf{14} 560--565.

\bibitem{GuiWang2012}
    Guillin, A., Wang, F.-Y. (2012).
    Degenerate Fokker-Planck Equations : Bismut Formula, Gradient Estimate and Harnack Inequality.
    \textit{J. Differential Equations},
    \textbf{253} 20--40.

\bibitem{Hsu2002}
    Hsu, E.~P. (2002).
    Stochastic Analysis on Manifolds.
    \textit{American Mathematical Society.}

\bibitem{Ouyang2009a}
    Ouyang, S.-X. 2009.
     Harnack Inequalities and Applications for Stochastic Equations, Ph.D thesis, Bielefeld University, 2009. Available on
\url{http://bieson.ub.uni-bielefeld.de/volltexte/2009/1463/pdf/ouyang.pdf}.

\bibitem{Ouyang2009b}
    Ouyang, S.-X. 2009.
    Non-time-homogeneous Generalized Mehler Semigroups and Applications,
    arXiv:1009.5314.

\bibitem{OuyangRW2012}
    Ouyang, S.-X. and R\"ockner, M. and Wang, F.-Y. 2012.
    Harnack inequalities and applications for Ornstein-Uhlenbeck semigroups with jump,
    \textit{Potential Anal.}
    \textbf{36}: 301--315.

\bibitem{RoLiu2010}
    R\"{o}ckner, M., Liu, W. (2010).
    SPDE in Hilbert space with locally monotone coefficients.
    \textit{J. Funct. Anal.}
    \textbf{259}  2902--2922.

\bibitem{ShaoWY2012}
    Shao, J., Wang, F.-Y. and Yuan, C. (2012).
    Harnack Inequalities for Stochastic (Functional) Differential Equations with Non-Lipschitzian Coefficients.
    \textit{arXiv:1208.5094}.

\bibitem{Wang11}
    Wang, F.-Y. (2011).
    Analysis on path spaces over Riemannian manifolds with boundary.
    \textit{Comm. Math. Sci.}
    \textbf{9} 1203--1212.

\bibitem{Wang2011}
    Wang, F.-Y. (2011).
    Harnack inequality for SDE with multiplicative noise and extension to Neumann semigroup on nonconvex manifolds.
    \textit{Ann. Probab.}
    \textbf{39} 1449--1467.

\bibitem{Wang2012}
    Wang, F.-Y. (2012).
    Integration by part formula and shift Harnack inequality for stochastic equation.
    \textit{arXiv:1203.4023}

\bibitem{WangWuXu2011}
    Wang, F.-Y., Wu, J.-L., Xu, L. (2011).
    Log-Harnack Inequality for Stochastic Burgers Equations and Applications.
    \textit{J. Math. Anal. Appl.}
    \textbf{384}  151--159.

\bibitem{WangXu2012}
    Wang, F.-Y., Xu, L.
    Derivative Formula and Applications for Hyperdissipative Stochastic Navier-Stokes/Burgers Equations.
    \textit{Infin. Dimens. Anal. Quantum Probab. Relat. Top.}(to appear).

\bibitem{WangZTS2010}
    Wang, F.-Y., Zhang, T.-S. (2010).
    Gradient estimates for stochastic evolution equations with non-Lipschitz coefficients.
    \textit{J. Math. Anal. Appl.}
    \textbf{365}  1--11.


\end{thebibliography}
\end{document}